\documentclass[12pt,reqno]{amsart}
\usepackage{amsfonts}
\usepackage{amssymb}
\usepackage{amscd}

\newtheorem{theorem}{Theorem}[section]

\newtheorem{proposition}[theorem]{Proposition}

\newtheorem{example}[theorem]{Example}
\numberwithin{equation}{section}
\begin{document}
\title[Lie groups]{Desingularizing compact Lie group actions}
\address{Department of Mathematics \\
Texas Christian University \\
Fort Worth, Texas 76129, USA}
\author{Ken Richardson}
\email{k.richardson@tcu.edu}

\begin{abstract}
This note surveys the well-known structure of $G$-manifolds and summarizes
parts of two papers that have not yet appeared: \cite{BKR}, joint with J. Br%
\"{u}ning and F. W. Kamber, and \cite{PrRi}, joint with I. Prokhorenkov. In
particular, from a given manifold on which a compact Lie group acts
smoothly, we construct a sequence of manifolds on which the same Lie group
acts, but with fewer levels of singular strata. Global analysis and
geometric results on the simpler manifolds may be translated to results on
the original manifold. Further, we show that by utilizing bundles over a $G$%
-manifold with singular strata, we may construct natural equivariant
transverse Dirac-type operators that have properties similar to Dirac
operators on closed manifolds.
\end{abstract}
\maketitle

\section{Manifolds and compact Lie group actions}

Suppose that a compact Lie group $G$ acts smoothly on a smooth, connected,
closed manifold $M$. We assume that the action is effective, meaning that no 
$g\in G$ fixes all of $M$. (Otherwise, replace $G$ with $G\diagup \left\{
g\in G:gx=x\text{ for all }x\in M\right\} =\emptyset $.) Choose a Riemannian
metric for which $G$ acts by isometries; average the pullbacks of any fixed
Riemannian metric over the group of diffeomorphisms to obtain such a metric.

Given such an action and $x\in M$, the \textbf{isotropy subgroup} $G_{x}<G$
is defined to be $\left\{ g\in G:gx=x\right\} $. The \textbf{orbit} $%
\mathcal{O}_{x}$ of a point $x$ is defined to be $\left\{ gx:g\in G\right\} $%
. Note that $G_{gx}=gG_{x}g^{-1}$, so the conjugacy class of the isotropy
subgroup of a point is fixed along an orbit. The conjugacy class of the
isotropy subgroups along an orbit is called the \textbf{orbit type}. On any
such $G$-manifold, there are a finite number of orbit types, and there is a
partial order on the set of orbit types. Given subgroups $H$ and $K$ of $G$
that occur as isotropy subgroups, we say that $\left[ H\right] \leq $ $\left[
K\right] $ if $H$ is conjugate to a subgroup of $K$, and we say $\left[ H%
\right] <$ $\left[ K\right] $ if $\left[ H\right] \leq $ $\left[ K\right] $
and $\left[ H\right] \neq $ $\left[ K\right] $. We may enumerate the
conjugacy classes of isotropy subgroups as $\left[ G_{0}\right] ,...,\left[
G_{r-1}\right] $ such that $\left[ G_{i}\right] \leq \left[ G_{j}\right] $
if and only if $i\leq j$. It is well-known that the union $M_{0}$ of the 
\textbf{principal orbits} (those with type $\left[ G_{0}\right] $) form an
open dense subset $M_{0}$ of the manifold $M$, and the other orbits are
called \textbf{singular}. As a consequence, every isotropy subgroup $H$
satisfies $\left[ G_{0}\right] \leq \left[ H\right] $; $M_{0}$ is called the 
\textbf{principal stratum}. Let $M_{j}$ denote the set of points of $M$ of
orbit type $\left[ G_{j}\right] $ for each $j$; the set $M_{j}$ is called
the \textbf{stratum} corresponding to $\left[ G_{j}\right] $. A stratum $%
M_{j}$ is called a \textbf{most singular stratum} if there does not exist a
stratum $M_{k}$ such that $\left[ G_{j}\right] <\left[ G_{k}\right] $. It is
known that each stratum is a $G$-invariant submanifold of $M $, and in fact
a most singular stratum is a closed (but not necessarily connected)
submanifold. Also, for each $j$, the submanifold $M_{\geq
j}:=\bigcup\limits_{\left[ G_{k}\right] \geq \left[ G_{j}\right] }M_{k}$ is
a closed, $G$-invariant submanifold.

Consider the following simple examples.

\begin{example}
$G=\mathbb{Z}_{2}\times \mathbb{Z}_{2}$ acts on $M=S^{2}\subset \mathbb{R}%
^{3}$ by 
\begin{equation*}
\left( 1,0\right) \left( x,y,z\right) =\left( x,-y,z\right) ,~\left(
0,1\right) \left( x,y,z\right) =\left( x,y,-z\right) .
\end{equation*}
The action is isometric for the standard metric on the sphere. The isotropy
subgroup of each point in $M_{0}=\left\{ \left( x,y,z\right) \in S^{2}:z\neq
0,y\neq 0\right\} $ is $G_{0}=\left\{ \left( 0,0\right) \right\} $, points
of the set $M_{2}=\left\{ \left( x,y,z\right) \in S^{2}:z=0,y\neq 0\right\} $
have isotropy subgroup $G_{2}=\left\{ \left( 0,0\right) ,\left( 0,1\right)
\right\} $, points in $M_{1}=\left\{ \left( x,y,z\right) \in S^{2}:y=0,z\neq
0\right\} $ have isotropy subgroup $G_{1}=\left\{ \left( 0,0\right) ,\left(
1,0\right) \right\} $, and points of $M_{3}=\left\{ \left( 1,0,0\right)
,\left( -1,0,0\right) \right\} $ have isotropy subgroup $G_{3}=\mathbb{Z}%
_{2}\times \mathbb{Z}_{2}$. Here, $M_{0}$ is the principal stratum, and $%
M_{3}$ is the only most singular stratum. The set $M_{\geq 1}$ is the circle 
$\left\{ \left( x,y,z\right) \in S^{2}:y=0\right\} $. Note that $\left[ G_{0}%
\right] <\left[ G_{1}\right] <\left[ G_{3}\right] $ and $\left[ G_{0}\right]
<\left[ G_{2}\right] <\left[ G_{3}\right] $ but $\left[ G_{1}\right] $ and $%
\left[ G_{2}\right] $ are not comparable.
\end{example}

\begin{example}
$G=S^{1}$ acts on $M=S^{2}\subset \mathbb{R}^{3}$ by 
\begin{equation*}
e^{i\theta }\left( x,y,z\right) =\left( x\cos \left( \theta \right) -y\sin
\left( \theta \right) ,x\sin \left( \theta \right) +y\cos \left( \theta
\right) ,z\right) ,
\end{equation*}
i.e. rotations around the $z$-axis. This action is isometric for the
standard metric. The isotropy subgroup at each point of $M_{0}=\left\{
\left( x,y,z\right) \in S^{2}:\left\vert z\right\vert <1\right\} $ is $%
G_{0}=\left\{ 1\right\} $, and the isotropy subgroup for each point of $%
M_{1}=\left\{ \left( x,y,z\right) \in S^{2}:z=1\text{ or }z=-1\right\} $ is $%
G_{1}=S^{1}$.
\end{example}

\begin{example}
\label{torusAction}$G=S^{1}\times S^{1}$ acts on $M=S^{3}\subset \mathbb{R}%
^{4}$ by%
\begin{multline*}
\left( e^{i\theta },e^{i\alpha }\right) \left( x,y,z,w\right)  \\
=(~x\cos \left( \theta \right) -y\sin \left( \theta \right) ,x\sin \left(
\theta \right) +y\cos \left( \theta \right) ,z\cos \left( \alpha \right)
-w\sin \left( \alpha \right) ,z\sin \left( \alpha \right) +w\cos \left(
\alpha \right) ~).
\end{multline*}
The isotropy subgroup of each point of 
\begin{equation*}
M_{0}=\left\{ \left( x,y,z,w\right) \in
S^{3}:x^{2}+y^{2}>0,z^{2}+w^{2}>0\right\} 
\end{equation*}
is $G_{0}=\left\{ \left( 1,1\right) \right\} $, the isotropy subgroup for
each point of $M_{1}=\left\{ \left( x,y,z,w\right) \in S^{3}:z=w=0\right\} $
is $G_{1}=\left\{ 1\right\} \times S^{1}$, and the isotropy subroup of $%
M_{2}=\left\{ \left( x,y,z,w\right) \in S^{3}:x=y=0\right\} $ is $%
G_{2}=S^{1}\times \left\{ 1\right\} $. Note that the diagonal action of $%
\left( e^{i\theta },e^{i\theta }\right) $ with $e^{i\theta }\in S^{1}$ gives
the Hopf fibration. In this example, the torus action on $S^{3}$ yields
three strata, with $M_{0}$ being the principal stratum and $M_{1}$ and $M_{2}
$ each being a most singular stratum.
\end{example}

\section{Desingularization construction}

With notation as in the previous section, we will construct a new $G$%
-manifold $N$ that has a single stratum (of type $\left[ G_{0}\right] $) and
that is a branched cover of $M$, branched over the singular strata. A
distinguished fundamental domain of $M_{0}$ in $N$ is called the
desingularization of $M$ and is denoted $\widetilde{M}$. The significance of
this construction is that it appears in the equivariant index theorem in 
\cite{BKR}, and the analysis of transversally elliptic operators on $M$ may
be replaced by analysis on $\widetilde{M}$, which is much easier to
understand.

A sequence of constructions is used to construct $N$ and $\widetilde{M}%
\subset N$. If $M_{j}$ is a most singular stratum, let $T_{\varepsilon
}\left( M_{j}\right) $ denote an open tubular neighborhood of $M_{j}$ of
radius $\varepsilon >0$. If $\varepsilon $ is sufficiently small, then all
orbits in $T_{\varepsilon }\left( M_{j}\right) \setminus M_{j}$ are of type $%
\left[ G_{k}\right] $, where $\left[ G_{k}\right] <\left[ G_{j}\right] $.
Let 
\begin{equation*}
N^{1}=\left( M\setminus T_{\varepsilon }\left( M_{j}\right) \right) \cup
_{\partial T_{\varepsilon }\left( M_{j}\right) }\left( M\setminus
T_{\varepsilon }\left( M_{j}\right) \right)
\end{equation*}%
be the manifold constructed by gluing two copies of $\left( M\setminus
T_{\varepsilon }\left( M_{j}\right) \right) $ smoothly along the boundary.
Since the $T_{\varepsilon }\left( M_{j}\right) $ is saturated (a union of $G$%
-orbits), the $G$-action lifts to $N^{1}$. Note that the strata of the $G$%
-action on $N^{1}$ correspond to strata in $M\setminus T_{\varepsilon
}\left( M_{j}\right) $. If $M_{k}\cap \left( M\setminus T_{\varepsilon
}\left( M_{j}\right) \right) $ is nontrivial, then the stratum corresponding
to isotropy type $\left[ G_{k}\right] $ on $N^{1}$ is 
\begin{equation*}
N_{k}^{1}=\left( M_{k}\cap \left( M\setminus T_{\varepsilon }\left(
M_{j}\right) \right) \right) \cup _{\left( M_{k}\cap \partial T_{\varepsilon
}\left( M_{j}\right) \right) }\left( M_{k}\cap \left( M\setminus
T_{\varepsilon }\left( M_{j}\right) \right) \right) .
\end{equation*}%
Thus, $N^{1}$ is a $G$-manifold with one fewer stratum than $M$, and $%
M\setminus M_{j}$ is diffeomorphic to one copy of $\left( M\setminus
T_{\varepsilon }\left( M_{j}\right) \right) $, denoted $\widetilde{M}^{1}$
in $N^{1}$. In fact, $N^{1}$ is a branched double cover of $M$, branched
over $M_{j}$. If $N^{1}$ has one orbit type, then we set $N=N^{1}$ and $%
\widetilde{M}=\widetilde{M}^{1}$. If $N^{1}$ has more than one orbit type,
we repeat the process with the $G$-manifold $N^{1}$ to produce a new $G$%
-manifold $N^{2}$ with two fewer orbit types than $M$ and that is a $4$-fold
branched cover of $M$. Again, $\widetilde{M}^{2}$ is a fundamental domain of 
$\widetilde{M}^{1}\setminus \left\{ \text{a most singular stratum}\right\} $%
, which is a fundamental domain of $M$ with two strata removed. We continue
until $N=N^{r-1}$ is a $G$-manifold with all orbits of type $\left[ G_{0}%
\right] $ and is a $2^{r-1}$-fold branched cover of $M$, branched over $%
M\setminus M_{0}$. We set $\widetilde{M}=\widetilde{M}^{r-1}$, which is a
fundamental domain of $M_{0}$ in $N$.

As mentioned earlier, if $M$ is equipped with a $G$-equivariant,
transversally elliptic differential operator on sections of an equivariant
vector bundle over $M$, then this data may be pulled back to the
desingularization $\widetilde{M}$. Given the bundle and operator over $N^{j}$%
, simply form the invertible double of the operator on $N^{j+1}$, which is
the double of the manifold with boundary $N^{j}\setminus T_{\varepsilon
}\left( \Sigma \right) $, where $\Sigma $ is a most singular stratum on $%
N^{j}$.

Further, one may independently desingularize $M_{\geq j}$, since this
submanifold is itself a closed $G$-manifold. If $M_{\geq j}$ has more than
one connected component, we may desingularize all components simultaneously.
Note that the isotropy type of all points of $\widetilde{M_{\geq j}}$ is $%
\left[ G_{j}\right] $, and the $\widetilde{M_{\geq j}}\diagup G$ is a smooth
(open) manifold.

The desingularizations $\widetilde{M}$ and $\widetilde{M_{\geq j}}$ are the
regions of integration present in the following equivariant index formula in 
\cite{BKR}.

\begin{theorem}
(Equivariant Index Theorem, in \cite{BKR})\label{EquivariantIndexTheorem}%
\begin{eqnarray*}
\mathrm{ind}^{\rho }\left( D\right) &=&\int_{\widetilde{M}\diagup
G}A_{0}^{\rho }\left( x\right) ~\widetilde{\left\vert dx\right\vert }~ \\
&&+\sum_{j,a,b}C_{jab}\int_{\widetilde{M_{\geq j}}\diagup G}\left( -\eta
\left( D_{j}^{S+,\sigma _{a}}\right) +h\left( D_{j}^{S+,\sigma _{a}}\right)
\right) A_{j,\sigma _{b}^{\ast }}^{\rho _{0}}\left( x\right) ~\widetilde{%
\left\vert dx\right\vert }.
\end{eqnarray*}
\end{theorem}

In this formula, $D$ is a $G$-equivariant, transversally elliptic operator
which can be written in a tubular neighborhood of each $M_{j}$ as the
product 
\begin{equation*}
D=\left\{ Z_{j}\left( \nabla _{\partial _{r}}^{E}+\frac{1}{r}%
D_{j}^{S}\right) \right\} \ast D^{M_{\geq j}},
\end{equation*}%
where $r$ is the distance from $M_{\geq j}$, where $Z_{j}$ is a local bundle
isomorphism, the map $D_{j}^{S}$ is a family of purely first order operators
that differentiates in the unit normal bundle directions tangent to $%
S_{x}M_{\geq j}$, and $D^{M_{\geq j}}\ $is a global transversally elliptic, $%
G$-equivariant, first order operator on the stratum $M_{\geq j}$. Many
important examples of equivariant, transversally elliptic differential
operators satisfy the condition above. In the formula, the forms $%
A_{0}^{\rho }$ and $A_{j,\sigma _{b}^{\ast }}^{\rho _{0}}$ are Atiyah-Singer
integrands for differential operators on the corresponding desingularized
strata. The number $\eta \left( D_{j}^{S+,\sigma _{a}}\right) $ corresponds
to a $G_{j}$-equivariant eta invariant, and $h\left( D_{j}^{S+,\sigma
_{a}}\right) $ is the dimension of the $\sigma _{a}$-part of the kernel of $%
D_{j}^{S+}$, which is a $G_{j}$ representation space. These numbers are
actually topological invariants, a little surprising if one is familiar with
eta invariants of operators on closed manifolds; the invariance comes from
the fact that the polar decomposition yields integral eigenvalues for $%
D_{j}^{S}$, and thus they are constant along connected components of $M_{j}$%
. The constants $C_{jab}$ depend on the representation theory --- the
induced representation of $G_{j}$ on a $\rho $-representation space and on
Clebsch-Gordan coefficients of tensor product representations.

\section{Natural equivariant Dirac operators\label{naturalSection}}

Here another approach is used to treat difficult transverse analytic
problems in a less singular setting, and the content of this section is in
the paper \cite{PrRi}. Given a connected, complete $G$-manifold, the action of 
$g\in G$ on $M$ induces an action of $dg$ on $TM$, which in turn induces an
action of $G$ on the principal $O\left( n\right) $-bundle $F_{O}\overset{p}{%
\rightarrow }M$ of orthonormal frames over $M$. It turns out that when $G$
is effective, the isotropy subgroups on $F_{O}$ are all trivial. In any
case, the $G$ orbits on $F_{O}$ are diffeomorphic and are the fibers
(leaves) of a Riemannian fiber bundle in a natural Sasakian metric on $F_{O}$%
. The quotient $F_{O}\overset{\pi }{\rightarrow }F_{O}\diagup G$ is a
Riemannian submersion of compact $O\left( n\right) $-manifolds. The metric
on $F_{O}$ is bundle-like for the Riemannian foliation $\left( F_{O},%
\mathcal{F}\right) $ by $G$-orbits.

Let $E\rightarrow F_{O}$ be a Hermitian vector bundle that is equivariant
with respect to the $G\times O\left( n\right) $ action. Let $\rho
:G\rightarrow U\left( V_{\rho }\right) $ and $\sigma :O\left( n\right)
\rightarrow U\left( W_{\sigma }\right) $ be irreducible unitary
representations. We define the bundle $\mathcal{E}^{\sigma }\rightarrow M$
by 
\begin{equation*}
\mathcal{E}_{x}^{\sigma }=\Gamma \left( p^{-1}\left( x\right) ,E\right)
^{\sigma },
\end{equation*}%
where the superscript $\sigma $ refers to the type $\sigma $ part of the $%
O\left( n\right) $-representation space $\Gamma \left( p^{-1}\left( x\right)
,E\right) $. The bundle $\mathcal{E}^{\sigma }$ is a Hermitian $G$-vector
bundle of finite rank over $M$ that comes with a natural metric.

Similarly, we define the bundle $\mathcal{T}^{\rho }\rightarrow F_{O}\diagup
G$ by%
\begin{equation*}
\mathcal{T}_{y}^{\rho }=\Gamma \left( \pi ^{-1}\left( y\right) ,E\right)
^{\rho },
\end{equation*}%
and $\mathcal{T}^{\rho }\rightarrow F_{O}\diagup G$ is a Hermitian $O\left(
n\right) $-equivariant bundle of finite rank, with a natural metric.

\begin{theorem}
\label{IsomorphismsOfSectionsTheorem}(in \cite{PrRi}) For any irreducible
representation $\rho :G\rightarrow U\left( V_{\rho }\right) $ and
irreducible representation $\sigma :O\left( n\right) \rightarrow U\left(
W_{\sigma }\right) $, there is an explicit isomorphism $\Gamma \left( M,%
\mathcal{E}^{\sigma }\right) ^{\rho }\rightarrow \Gamma \left( F_{O}\diagup
G,\mathcal{T}^{\rho }\right) ^{\sigma }$ that extends to an $L^{2}$-isometry.
\end{theorem}

Next, let $E\rightarrow F_{O}$ be a Hermitian vector bundle of $\mathbb{C}%
\mathrm{l}\left( N\mathcal{F}\right) $ modules that is equivariant with
respect to the $G\times O\left( n\right) $ action. As explained in \cite%
{DGKY}, \cite{GlK}, \cite{PrRi}, \cite{BKR}, \cite{HabRi}, there is natural
construction of what is know as the \textbf{basic Dirac operator} in this
situation. We have the transversal Dirac operator $D_{\mathrm{tr}}$ defined
by the composition 
\begin{equation*}
\Gamma \left( F_{O},E\right) \overset{\nabla }{\rightarrow }\Gamma \left(
F_{O},T^{\ast }F_{O}\otimes E\right) \overset{\mathrm{proj}}{\rightarrow }%
\Gamma \left( F_{O},N^{\ast }\mathcal{F}\otimes E\right) \overset{c}{%
\rightarrow }\Gamma \left( F_{O},E\right) .
\end{equation*}%
The basic Dirac operator 
\begin{equation*}
D_{N\mathcal{F}}=\frac{1}{2}\left( D_{\mathrm{tr}}+D_{\mathrm{tr}}^{\ast
}\right) =D_{\mathrm{tr}}-\frac{1}{2}c\left( H\right)
\end{equation*}%
is a essentially self-adjoint $G\times O\left( n\right) $-equivariant
operator, where $H$ is the mean curvature vector field of the $G$-orbits in $%
F_{O}$.

From $D_{N\mathcal{F}}$ we now construct equivariant differential operators
on $M$ and $F_{O}\diagup G$, denoted%
\begin{equation*}
D_{M}^{\sigma }:\Gamma \left( M,\mathcal{E}^{\sigma }\right) \rightarrow
\Gamma \left( M,\mathcal{E}^{\sigma }\right)
\end{equation*}%
and%
\begin{equation*}
D_{F_{O}\diagup G}^{\rho }:\Gamma \left( F_{O}\diagup G,\mathcal{T}^{\rho
}\right) \rightarrow \Gamma \left( F_{O}\diagup G,\mathcal{T}^{\rho }\right)
,
\end{equation*}%
defined in the natural way. For an irreducible representation $\alpha
:G\rightarrow U\left( V_{\alpha }\right) $, let 
\begin{equation*}
\left( D_{M}^{\sigma }\right) ^{\alpha }:\Gamma \left( M,\mathcal{E}^{\sigma
}\right) ^{\alpha }\rightarrow \Gamma \left( M,\mathcal{E}^{\sigma }\right)
^{\alpha }
\end{equation*}%
be the restriction of $D_{M}^{\sigma }$ to sections of $G$-representation
type $\left[ \alpha \right] $. Similarly, for an irreducible representation $%
\beta :G\rightarrow U\left( W_{\beta }\right) $, let 
\begin{equation*}
\left( D_{F_{O}\diagup G}^{\rho }\right) ^{\beta }:\Gamma \left(
F_{O}\diagup G,\mathcal{T}^{\rho }\right) ^{\beta }\rightarrow \Gamma \left(
F_{O}\diagup G,\mathcal{T}^{\rho }\right) ^{\beta }
\end{equation*}%
be the restriction of $D_{F_{O}\diagup G}^{\rho }$ to sections of $O\left(
n\right) $-representation type $\left[ \beta \right] $. The proposition
below follows from Theorem \ref{IsomorphismsOfSectionsTheorem}.

\begin{proposition}
The operator $D_{M}^{\sigma }$ is transversally elliptic and $G$%
-equivariant, and $D_{F_{O}\diagup G}^{\rho }$ is elliptic and $O\left(
n\right) $-equivariant, and the closures of these operators are
self-adjoint. The operators $\left( D_{M}^{\sigma }\right) ^{\rho }$ and $%
\left( D_{F_{O}\diagup G}^{\rho }\right) ^{\sigma }$ have identical discrete
spectrum, and the corresponding eigenspaces are conjugate via Hilbert space
isomorphisms.
\end{proposition}

Thus, questions about the transversally elliptic operator $D_{M}^{\sigma }$
are reduced to questions about the elliptic operators $D_{F_{O}\diagup
G}^{\rho }$ for each irreducible $\rho :G\rightarrow U\left( V_{\rho
}\right) $. Further, it turns out that the operators $D_{M}^{\sigma }$ play
the same role for equivariant analysis as the standard Dirac operators do in
the index theory and analysis of elliptic operators on closed manifolds.

\section{Further Comments}

The main use of these results is that in both cases, a quite difficult
problem of analyzing a transversally elliptic operator (with potentially
infinite dimensional eigenspaces) is reduced to an elliptic problem or set
of elliptic problems, which are much more tractable. For example, the
Atiyah-Segal Theorem (\cite{ASe}) was the first version of an equivariant
index theorem, and it appeared in 1968. However, the appropriate
generalization to transversally elliptic operators appeared only in 1996 and
was due to Berline and Vergne (\cite{Be-V1},\cite{Be-V2}). In a sense,
Theorem \ref{EquivariantIndexTheorem} is a Fourier transform version of the
Atiyah-Segal and Berline-Vergne results, giving a formula for the Fourier
coefficients of the character instead of the value of the character at a
particular $g\in G$. Further, Theorem \ref{EquivariantIndexTheorem} gives a
method of computing eta invariants of Dirac-type operators on quotients of
spheres by compact group actions; this was known previously for finite group
actions only.

The new \textquotedblleft transversal Dirac operators\textquotedblright\ on $%
G$-manifolds constructed in Section \ref{naturalSection} and in \cite{PrRi}
should be explored further, and in particular future investigations should
lead to generalizations of Dirac operator results to the transversally
elliptic setting.

\end{document}